\documentclass[11pt,french,english]{amsart}
\usepackage[T1]{fontenc}
\usepackage[latin1]{inputenc}
\usepackage{latexsym}
\usepackage{textcomp}
\usepackage{amsthm}
\usepackage{relsize}
\usepackage{amssymb}

\makeatletter

\DeclareRobustCommand{\lyxmathsym}[1]{\ifmmode\begingroup\def\b@ld{bold}
  \def\rmorbf##1{\ifx\math@version\b@ld\textbf{##1}\else\textrm{##1}\fi}
  \mathchoice{\hbox{\rmorbf{#1}}}{\hbox{\rmorbf{#1}}}
  {\hbox{\smaller[2]\rmorbf{#1}}}{\hbox{\smaller[3]\rmorbf{#1}}}
  \endgroup\else#1\fi}

\providecommand{\tabularnewline}{\\}

\numberwithin{equation}{section} %% Comment out for sequentially-numbered
\numberwithin{figure}{section} %% Comment out for sequentially-numbered
\theoremstyle{plain}
\theoremstyle{plain}
\newtheorem{thm}{Theorem}
  \theoremstyle{definition}
  \newtheorem{condition}[thm]{Hypothesis}
  \theoremstyle{plain}
  \newtheorem{lem}[thm]{Lemma}
  \theoremstyle{definition}
  \newtheorem{defn}[thm]{Definition}
  \theoremstyle{plain}
  \newtheorem{cor}[thm]{Corollary}
  \theoremstyle{remark}
  \newtheorem{rem}[thm]{Remark}
  \theoremstyle{plain}
  \newtheorem{prop}[thm]{Proposition}

\AtBeginDocument{
  
}

\makeatother

\usepackage{babel}
\addto\extrasfrench{\providecommand{\fg}{\ifdim\lastskip>\z@\unskip\fi~\frqq}}

\begin{document}
\newcommand{\Alb}{{\rm Alb}}

\newcommand{\Jac}{{\rm Jac}}

\newcommand{\Hom}{{\rm Hom}}

\newcommand{\End}{{\rm End}}

\newcommand{\Aut}{{\rm Aut}}

\newcommand{\NS}{{\rm NS}}

\title{Genus $2$ curve configurations on Fano surfaces}

\author{Xavier Roulleau}
\begin{abstract}
We study the configurations of genus $2$ curves on the Fano surfaces
of cubic threefolds. We establish a link between some involutive automorphisms
acting on such a surface $S$ and genus $2$ curves on $S$. We give
a partial classification of the Fano surfaces according to the automorphism
group generated by these involutions and determine the configurations
of their genus $2$ curves. We study the Fano surface of the Klein
cubic threefold for which the $55$ genus $2$ curves generate a rank
$25=h^{1,1}$ index $2$ subgroup of the Néron-Severi group.
\end{abstract}
\maketitle


\section{Introduction.}

Let $S$ be a smooth surface which verifies the following Hypothesis:
\begin{condition}
\label{l'Hypoth=0000E8se}The variety $S$ is a smooth complex surface
of general type. The cotangent sheaf $\Omega_{S}$ of $S$ is generated
by its space $H^{0}(\Omega_{S})$ of global sections and the irregularity
$q=\dim H^{0}(\Omega_{S})$ satisfies $q>3$.
\end{condition}
Let $T_{S}$ be the tangent sheaf, $\pi:\mathbb{P}(T_{S})\rightarrow S$
be the projection and let \[
\psi:\mathbb{P}(T_{S})\rightarrow\mathbb{P}(H^{o}(\Omega_{S})^{*})=\mathbb{P}^{q-1}\]
be the cotangent map of $S$ defined by : $\pi_{*}\psi^{*}\mathcal{O}_{\mathbb{P}^{q-1}}(1)=\Omega_{S}$.\\
In \cite{Roulleau} a curve $C\hookrightarrow S$ is called non-ample
if there is a section \[
t:C\rightarrow\mathbb{P}(T_{S})\]
 such that $\psi(t(C))$ is a point. The cotangent sheaf of $S$ is
ample (i.e. $\psi^{*}\mathcal{O}_{\mathbb{P}^{q-1}}(1)$ is ample)
if and only if $S$ does not contain non-ample curves \cite{Roulleau}.
In other words, the non-ample curves are the obstruction for the cotangent
sheaf to be ample. A natural example of a non-ample curve is given
by a smooth curve $C\hookrightarrow S$ of genus $1$, where a section
$t:C\rightarrow\mathbb{P}(T_{S})$ such that $\psi(t(C))$ is a point
is given by the natural quotient: \[
\Omega_{S}\otimes\mathcal{O}_{C}\rightarrow\Omega_{C}.\]

A step further after the study of ampleness of $\psi^{*}\mathcal{O}_{\mathbb{P}^{q-1}}(1)$
is the study of the very-ampleness. In projective geometry, the simplest
object after points are lines. An obstruction for invertible sheaf
$\psi^{*}\mathcal{O}_{\mathbb{P}^{q-1}}(1)$ to be very ample is that
the surface contains a curve $C$ for which there is a section $t:C\rightarrow\mathbb{P}(T_{S})$
such that $\psi(t(C))$ is a line. We say that such a curve $C$ satisfies
property $(*)$. A natural example of such a curve is given by a smooth
curve $C\hookrightarrow S$ genus $2$, where the section $t:C\rightarrow\mathbb{P}(T_{S})$
is given by the natural quotient: \[
\Omega_{S}\otimes\mathcal{O}_{C}\rightarrow\Omega_{C}.\]

In the present paper, we classify the curves that satisfy property
$(*)$ on Fano surfaces and we focus our attention on the genus $2$
curves. \\
By definition, a Fano surface parametrizes the lines of a smooth
cubic threefold $F\hookrightarrow\mathbb{P}^{4}$. This scheme is
a surface $S$ that verifies Hypothesis \ref{l'Hypoth=0000E8se} and
has irregularity $q=5$. By the Tangent Bundle Theorem \cite{Clemens}
12.37, the image $F'$ of the cotangent map $\psi:\mathbb{P}(T_{S})\rightarrow\mathbb{P}(H^{0}(\Omega_{S})^{*})$
of $S$ is a hypersurface of $\mathbb{P}(H^{0}(\Omega_{S})^{*})\simeq\mathbb{P}^{4}$
that is isomorphic to the original cubic $F$. Moreover, when we identify
$F$ and $F'$, the triple $(\mathbb{P}(T_{S}),\pi,\psi)$ is the
universal family of lines on $F$.\\
By a generic point of the cubic threefold $F$ goes $6$ lines.
For $s$ a point of $S$ the Fano surface of lines on $F$, we denote
by $L_{s}$ the corresponding line on $F$ and we denote by $C_{s}$
theincidence divisor :\[
C_{s}=\overline{\{t/t\not=s,\, L_{t}\,\mbox{cuts}\, L_{s}\}}.\]
This divisor is smooth and has genus $11$ for $s$ generic.

Let $D\hookrightarrow S$ be a curve such that there is a section
$D\rightarrow\mathbb{P}(T_{S})$ mapped onto a line by $\psi$ (property
$(*)$). There exists a point $t$ of $S$ such that this line is
$L_{t}$. As all the lines of $\psi_{*}\pi^{*}D$ goes through $L_{t}$,
there exists a residual divisor $R$ such that: \[
C_{t}=D+R.\]

\begin{thm}
\label{thm:A-curve-D}A curve $D$ on $S$ that verifies property
$(*)$ is a non-genus $1$ irreducible component of an incidence divisor
$C_{t}$.

If an incidence divisor $C_{t}$ is not irreducible, then it can split
only in the following way:\\
a) $C_{t}=E+R$ where $E$ is an elliptic curve, $R$ verifies
property $(*)$ and has genus $7$, it is an irreducible fiber of
a fibration of $S$ and $RE=4$.\\
b) $C_{t}=E+R+E'$ where $E,E'$ are two smooth curves of genus
$1$. The curve $R$ verifies property $(*)$ and has arithmetical
genus $4$.\\
c) $C_{t}=D+R$ where $D$ is a smooth genus $2$ curve, the curve
$R$ has arithmetical genus $4$, $D^{2}=-4$, $DR=6$ and $R^{2}=-3$.
The curves $D$ and $R$ satisfy property $(*)$.\\
d) If $D\not=D'\hookrightarrow S$ are two smooth curves of genus
$2$ on $S$, then $C_{s}D=2$ and :\[
DD'\in\{0,1,2\}.\]

\end{thm}
The involutive automorphisms of a Fano surface $S$ can be classified
into two type, say I and II. We prove in \cite{Roulleau} that there
is a bijection between the set of elliptic curves $E$ on $S$ and
the set of involutions $\sigma_{E}$ of type I, in such a way that:
\[
(\sigma_{E}\sigma_{E'})^{EE'}=1,\]
where $EE'$ is the intersection number of $E$ and $E'$. This formula
and the full classification of groups generated by involutions of
type I, enable us to determine all the configurations of genus $1$
curves on Fano surfaces. We give here an analogous result for automorphisms
of type II and genus $2$ curves:
\begin{thm}
To each involution $g$ of type II in $\Aut(S)$, there corresponds
a curve $D_{g}$ of genus $2$ on $S$. 

Let $G$ be one of the groups : $\mathbb{Z}/2\mathbb{Z},$ $\mathbb{D}_{n},\, n\in\{2,3,5,6\}$
(dihedral group of order $2n$), $\mathbb{A}_{5}$ (alternating group)
and $PSL_{2}(\mathbb{F}_{11})$. There exists a Fano surface $S$
with the following properties : \\
A) We can identify $G$ to a sub-group of $\Aut(S)$. By this identification,
each involution of $G$ has type II.\\
B) The number $N$ of genus $2$ curves on $S$ is as follows:\smallskip{}
\\
\begin{tabular}{|c|c|c|c|c|c|c|c|}
\hline 
Group $G$ & $\mathbb{Z}/2\mathbb{Z}$ & $\mathbb{D}_{2}$ & $\mathbb{D}_{3}$ & $\mathbb{D}_{5}$ & $\mathbb{D}_{6}$ & $\mathbb{A}_{5}$ & $PSL_{2}(\mathbb{F}_{11})$\tabularnewline
\hline 
$N$  & 1 & 3 & 3 & 5 & 7 & 15 & 55\tabularnewline
\hline
\end{tabular}\\
\\
These curves are smooth for $S$ generic.\\
C) The intersection number of the genus $2$ curves $D_{g},D_{h}$
is given by the formula:\[
D_{g}D_{h}=\left\{ \begin{array}{ccc}
-4 & if & g=h\\
0 & if & o(gh)=2\,\, or\,\,6\\
2 & if & o(gh)=3\\
1 & if & o(gh)=5\end{array}\right.\]
where o(f) is the order of the element $f$.\\
 D) For $G=\mathbb{D}_{3}$ (resp. $G=\mathbb{D}_{5}$), let $D$
be the sum of the $3$ (resp. $5$) genus $2$ curves $D_{g}$ ($g$
involution of $G$). The divisor $D$ is a fiber of a fibration of
$S$. \\
E) For $G=\mathbb{A}_{5}$, the $15$ genus $2$ curves generates
a sub-lattice $\Lambda$ of $\NS(S)$ of rank $15$, signature $(1,14)$
and discriminant $2^{24}3^{6}$. For $S$ generic, $\Lambda$ has
finite index inside $\NS(S)$. There exist an infinite number of such
surfaces with maximal Picard number $h^{1,1}=25$.\\
F) For $G=PSL_{2}(\mathbb{F}_{11})$, $S$ is the Fano surface
of the Klein cubic \[
x_{1}x_{5}^{2}+x_{5}x_{3}^{2}+x_{3}x_{4}^{2}+x_{4}x_{2}^{2}+x_{2}x_{1}^{2}=0.\]
The sublattice $\Lambda'$ of the Néron-Severi group $\NS(S)$ generated
by the $55$ smooth genus $2$ curves has rank $25=h^{1,1}$ and discriminant
$2^{2}11^{10}$. The group $\NS(S)$ is generated by $\Lambda'$ and
the class of an incidence divisor.
\end{thm}
We want to mention that Edge in \cite{Edge} was the first to give
a classification of cubics threefolds according to subgroups of $PSL_{2}(\mathbb{F}_{11})$,
in order to understand the geometry of a genus $26$ curve embedded
in $\mathbb{P}^{4}$, with automorphism group $PSL_{2}(\mathbb{F}_{11})$. 

After proving the above classification Theorem, we discuss on its
completeness. We finish this paper by a conjecture about the existence
of a surface of special type in $\mathbb{P}^{4}$ with a particular
configuration of $55$ $(-2)$-curves. The conjecture comes naturally
after studying the configuration of the $55$ genus $2$ curves on
the Fano surface of the Klein cubic.

\thanks{The author wish to gratefully acknowledge its host researcher Prof.
Miyaoka, and the JSPS organization for providing the support of this
work.}

\section{Genus $2$ curves on Fano surfaces.}

Let $S$ be the Fano surface of lines on a smooth cubic threefold
$F\hookrightarrow\mathbb{P}^{4}$. We denote by $L_{s}$ the line
on $F$ that corresponds to a point $s$ of $S$. Let $C_{s}$ be
the divisor that parametrizes the lines going through $L_{s}$.

Let $D\hookrightarrow S$ be a curve such that there is a section
$D\rightarrow\mathbb{P}(T_{S})$ mapped onto a line by the cotangent
map $\psi$ (property $(*)$). There exists a point $t$ of $S$ such
that this line is $L_{t}$. As all the lines of $\psi_{*}\pi^{*}D$
goes through $L_{t}$, there exists a residual divisor $R$ such that:
\[
C_{t}=D+R.\]
Let us prove Theorem \ref{thm:A-curve-D}. We first need to recall
some well-known facts on the Fano surfaces, the main reference used
here is  \cite{Bombieri}.\\
Let $s$ be a point of $S$ and $r$ be a generic point of $C_{s}$.
The plane $X$ that contains the lines $L_{s}$ and $L_{r}$ cuts
the cubic $F$ along the lines $L_{s},L_{r}$ and a third one denoted
by $L_{j_{s}(r)}$ such that :\[
XF=L_{s}+L_{r}+L_{j_{s}(r)}.\]
The rational map $j_{s}:C_{s}\rightarrow C_{s}$ extends to an automorphism
of $C_{s}$. Moreover the quotient $\Gamma_{s}$ of $C_{s}$ by $j_{s}$
parametrizes the planes $X$ that contain $L_{s}$ and that cut $F$
along three lines. The scheme $\Gamma_{s}$ is a plane quintic, and:
\begin{lem}
\label{lem:The-quotient}(\cite{Bombieri}, Lemma 2). The quotient
$j_{s}:C_{s}\rightarrow\Gamma_{s}$ is ramified over the singular
points of $\Gamma_{s}$. These singular points are ordinary double
points.
\end{lem}
We need also the following facts:
\begin{lem}
(\cite{Clemens}, Lemma 10.4, Proposition 10.3 formula 10.11). For
a point $s$ on $S$, the incidence divisor $C_{s}$ is ample, verifies
$C_{s}^{2}=5$ and $3C_{s}$ is numerically equivalent to the canonical
divisor of $S$.\\
(\cite{Roulleau}, Proposition 10, Theorem 13). Let $s$ be a point
of an elliptic curve $E\hookrightarrow S$. There exists an effective
divisor $R$ such that $C_{s}=E+R$. The curve $R$ is a fiber of
a fibration of $S$ onto $E$ and satisfies $RE=4$. Moreover, $E^{2}=-3$
and $C_{t}E=1$.
\end{lem}
Let us suppose that the quintic $\Gamma_{s}$ decomposes as $\Gamma_{s}=L+U$
where $L$ is a line. As $LU=4$, the component of $C_{s}$ over $L$
is smooth and ramified over $4$ points : it has genus $1$.\\
Conversely, suppose that $C_{t}=E+R$ where $E$ is an elliptic
curve. Then $ER=E(C_{s}-E)=4$ and $j_{s}$ restricted to $E$ is
a degree $2$ morphism ramified above $4$ points of $L$. Hence,
the image $L$ of $E$ by $j_{s}$ is a rational curve that cuts the
image of $R$ into $4$ points. Thus $L$ is a line.\\
Now, we proved in \cite{Roulleau} that an incidence divisor $C_{t}$
can have at most $2$ smooth genus $1$ irreducible components. Hence
there is at most $2$ lines that are components of the divisor $\Gamma_{s}$.
The case a) and b) of Theorem \ref{thm:A-curve-D} occur (see \cite{Roulleau}).

Suppose now that the quintic $\Gamma_{t}$ splits as $\Gamma_{t}=Q+T$
with $Q$ an irreducible quadric. Then the irreducible component $D$
of $C_{t}$ over $Q$ is branched over $6$ points of $Q$, thus this
is a smooth curve of genus $2$.

Let us recall that the Albanese map $S\rightarrow\Alb(S)$ is an embedding.
We consider $S$ as a subvariety of $\Alb(S)$. 
\begin{lem}
\label{lem: espace tangent ourbed}Let $D$ be a smooth curve of genus
$2$ on $S$, let $L_{t}\hookrightarrow F$ be the image by $\psi$
of the section obtained by $ $the quotient $\Omega_{S}\otimes\mathcal{O}_{D}\rightarrow\Omega_{D}$
and let $J(D)$ be the Jacobian of $D$. Consider the natural map
: $J(D)\hookrightarrow\Alb(S)$ . The tangent space \[
TJ(D)\hookrightarrow T\Alb(S)=H^{0}(\Omega_{S})^{*}\]
is the subjacent space to the line $L_{t}\hookrightarrow\mathbb{P}(H^{0}(\Omega_{S})^{*})=\mathbb{P}^{4}$,
i.e. $H^{0}(L_{t},\mathcal{O}_{L_{t}}(1))^{*}=TJ(D)$.\end{lem}
\begin{proof}
This is consequence of the definition of the cotangent map and the
fact that the section $D\rightarrow\mathbb{P}(T_{S})$ given by the
surjection \[
\Omega_{S}\otimes\mathcal{O}_{D}\rightarrow\Omega_{D}\]
is mapped onto the line $L_{t}$ by the cotangent map $\psi:\mathbb{P}(T_{S})\rightarrow F$.
\end{proof}
Let us suppose that the quintic $\Gamma_{t}$ decomposes as $\Gamma_{t}=Q+Q'+L$
with $Q,Q'$ irreducible quadrics and $L$ a line. Then there exist
$2$ smooth curves of genus $2$ $D,D'$ and an elliptic curve $E$
such that:\[
C_{t}=D+D'+E\]
Let $\vartheta:S\rightarrow\Alb(S)$ be an Albanese morphism of $S$.
By the previous Lemma, the $2$ curves $\vartheta(D),\vartheta(D')$
are contained on the same abelian surface $J(D)=J(D')\hookrightarrow\Alb(S)$.
Let $B$ be the quotient of $A$ by the smallest abelian variety containing
$J(D)$ and $E$ in $\Alb(S)$. The morphism $S\rightarrow\Alb(S)\rightarrow B$
contracts the ample divisor $C_{t}$ onto a point and its image generates
$B$ of dimension $>0$ : it is impossible.\\
 Thus if $\Gamma_{t}$ is reducible: $\Gamma_{t}=Q+T$ with $Q$
a smooth quadric, the other component $T$ is an irreducible cubic.

Let $R$ be the residual divisor of $D$ in the incidence divisor:
$C_{t}=D+R$. We have: \[
5=C_{t}^{2}=D^{2}+R^{2}+2DR.\]
 We know that $DR=6$ because the quadric $Q$ and the cubic $T$
such that $\Gamma_{t}=Q+T$ cut each others in $6$ points, thus :
\[
D^{2}+R^{2}=-7.\]
Moreover, $R^{2}=(C_{t}-D)^{2}=D^{2}-2C_{t}D+5$. Thus $2D^{2}-2C_{t}D+5=-7$
and \[
D^{2}-C_{t}D=-6\]
The divisor $3C_{t}$ is linearly equivalent to a canonical divisor
of $S$, thus $D^{2}+3C_{t}D=2$, and we obtain $C_{t}D=2$ and $D^{2}=-4$
and then $C_{t}R=3$, $R^{2}=-3$, $R$ has genus $4$.

Let $D'\not=D$  be a second smooth curve of genus $2$. As $DC_{s}=2$,
we obtain $D'(D+R)=2$. The intersection numbers $D'R$ and $DD'$
are positive, hence $0\leq DD'\leq2$.

\section{Link between genus $2$ curves and automorphisms, the Klein cubic.}

\subsection{Involutive automorphisms and genus $2$ curves on Fano surfaces.\protect \\
\label{sub:Involutive-automorphisms-and}}

Let $t$ be a point of $S$. We can suppose that the corresponding
line $L_{t}$ on the cubic is given by $x_{1}=x_{2}=x_{3}=0$, then
this cubic has equation :\[
\{C+2x_{4}Q_{1}+2x_{5}Q_{2}+x_{4}^{2}x_{1}+2x_{4}x_{5}\ell+x_{5}^{2}x_{3}=0\}\]
where $C,Q_{1},Q_{2},\ell$ are forms in the variables $x_{1},x_{2},x_{3}$.
Let $\mbox{\ensuremath{\Gamma}}_{t}$ be the scheme that parametrizes
the planes containing $L_{t}$ and such that their intersection with
$F$ is $3$ lines. This scheme $\Gamma_{t}$ is the quintic given
by the equation:\[
(x_{1}x_{3}-\ell^{2})C-Q_{1}^{2}x_{3}+2Q_{1}Q_{2}\ell-Q_{2}^{2}x_{1}=0\]
on the plane $\mathbb{P}(x_{1}:x_{2}:x_{3})$ (see \cite{Bombieri},
equation (6)). 
\begin{defn}
\label{def:An-automorphism-conjugated}An automorphism conjugated
to the involutive automorphism:\[
f:x\rightarrow(x_{1}:x_{2}:x_{3}:-x_{4}:-x_{5})\]
 in $PGL_{5}(\mathbb{C})$ is called a harmonic inversion of lines
and planes.
\end{defn}
The harmonic inversion $f$ acts on the cubic if and only if $Q_{1}=Q_{2}=0$.
In that case, a plane model of $\Gamma_{t}$ is \[
\Gamma_{t}=\{(x_{1}x_{3}-\ell^{2})C=0\}\]
and the cubic has equation:\[
F_{2}=\{C+x_{4}^{2}x_{1}+2x_{4}x_{5}\ell+x_{5}^{2}x_{3}=0\}.\]
If the conic $\mathcal{Q}=\{x_{1}x_{3}-\ell^{2}=0\}$ is smooth, the
divisor $C_{t}$, which is the double cover of $\Gamma{}_{t}$ branched
over the singularities of $\Gamma{}_{t}$, splits as follows:\[
C_{t}=D+R\]
with $D$ a smooth curve of genus $2$. If this conic $\mathcal{Q}$
is not smooth, then $C_{s}$ splits as follows: \[
C_{s}=E+E'+R\]
with $E$ and $E'$ two elliptic curves. Note that in the last case,
we can suppose $\ell=0$ and we see immediately that two other automorphisms
act on the cubic threefold. The divisor $E+E'$ has also genus $2$.
We proved:
\begin{cor}
\label{cor:To-each-order}To each harmonic inversion acting on the
cubic threefold $F$, there corresponds a curve of arithmetical genus
$2$ on the Fano surface of $F$. Such a curve is smooth or sum of
$2$ elliptic curves.
\end{cor}
Let $g$ be an harmonic inversion acting on $F$. It acts also on
$S$ and $H^{0}(\Omega_{S})$. 
\begin{lem}
\label{lem:The-trace-of}The trace of the action of $g$ on $H^{0}(\Omega_{S})$
is equal to $1$.\end{lem}
\begin{proof}
We can suppose that the cubic is :$ $\[
F_{2}=\{C+x_{4}^{2}x_{1}+2x_{4}x_{5}\ell+x_{5}^{2}x_{3}=0\}\]
and that $g:x\rightarrow(x_{1}:x_{2}:x_{3}:-x_{4}:-x_{5})$. By the
Tangent Bundle Theorem \cite{Clemens}, Theorem 12.37, we can consider
the homogeneous coordinates $x_{1},\dots,x_{5}$ as a basis of $H^{0}(\Omega_{S})$.
Thus, we see that the action of $g$ on $S$ has trace $1$ or $-1$.\\
The line $x_{1}=x_{2}=x_{3}=0$ lies in the cubic and correspond
to a fixed point $t$ of $f$. The action of $f$ on the tangent space
$TS_{t}=H^{0}(L_{t},\mathcal{O}(1))\hookrightarrow H^{0}(\Omega_{S})^{*}$
is thus the identity or the multiplication by $-1$. As $t$ is an
isolated fixed point, it is the multiplication by $-1$ and the action
of $g$ on $H^{0}(\Omega_{S})$ has trace $1$.\end{proof}
\begin{rem}
\label{rem: sur les theta caracteristiques}A) As a genus $2$ curve
has self-intersection number $-4$, there is a finite number of such
curves on a Fano surface.\\
B) It is equivalent to consider couples $(\Gamma,M)$ where $\Gamma$
is a plane quintic and $M$ is an odd theta characteristic of $\Gamma$
or to consider couples $(S,t)$ where $S$ a Fano surface and $t$
a point on $S$ (for these facts, see by example \cite{Casalaina},
Theorem 4.1).\\
 The moduli space of Fano surfaces with a smooth genus $2$ curve
is thus equal to the moduli of reducible plane quintics $\Gamma_{t}=Q+C$,
($Q$ smooth quadric, $C$ irreducible cubic), plus a theta characteristic.\\
The moduli space of $Q\simeq\mathbb{P}^{1}$ plus $6$ points $p_{1},\dots,p_{6}$
has dimension $3$. We can suppose that $ $$Q$ is the quadric $\{x^{2}+y^{2}+z^{2}=0\}$
inside $\mathbb{P}^{2}$. The cubics through $p_{1},\dots,p_{6}$
form a $3$ dimensional linear system. Thus the moduli space of Fano
surfaces that contains a genus $2$ curve is $6$ dimensional.\\
C) Given a reducible quintic curve $Q+C$ (with only simple singularities
and $Q$ smooth of degree $2$, $C$ irreducible), there are $32$
choices (\cite{Harris}, Corollary 2.7) of odd theta characteristics
giving non-isomorphic Fano surfaces containing a point $t$ such that
$\Gamma_{t}\simeq Q+C$. For only one of these Fano surfaces, the
genus $2$ curve corresponds to an order $2$ automorphism, namely
the Fano surface of the cubic: \[
F_{2}=\{C+x_{4}^{2}x_{1}+2x_{2}x_{4}x_{5}+x_{5}^{2}x_{3}=0\}.\]
D) As a genus $2$ curve is a cover of $\mathbb{P}^{1}$ branched
over $6$ points, any genus $2$ curve can be embedded inside a Fano
surface, in $\infty^{3}$ ways.\\
E) We have not found a geometric interpretation of the action of
an harmonic inversion on a Fano surface. That explains perhaps why
we do not know if our classification of group generated by harmonic
inversions acting on smooth cubic threefolds is complete or not.
\end{rem}

\subsection{Partial classification of configurations of genus $2$ curves, the
Fano surface of the Klein cubic threefold.\protect \\
}

Let $G$ be a group generated by order $2$ elements acting (faithfully)
on a $5$ dimensional space $V$ such that the trace of an order $2$
element is equal to $1$ (as in Lemma \ref{lem:The-trace-of}). \\
We say that a Fano surface $S$ (resp. the cubic threefold $F$
corresponding to $S$) has type $G$ if $G$ acts on $S$ and the
representation of $G$ on $H^{0}(\Omega_{S})$ is isomorphic to the
one on $V$. 

For the groups $\mathbb{D}_{n},\, n\in\{2,3,5,6\}$, $\mathbb{A}_{5}$
and $PSL_{2}(\mathbb{F}_{11})$, we take the unique $5$ dimensional
representation such that the elements have the following trace according
to their order:

\begin{tabular}{|c|c|c|c|c|c|}
\hline 
order & $2$ & $3$ & $5$ & $6$ & $11$\tabularnewline
\hline 
trace & 1 & -1 & 0 & 1 & $1/2(\lyxmathsym{\textminus}1+i\sqrt{11})$\tabularnewline
\hline
\end{tabular}

The aim of this paragraph is to prove the two following theorems:
\begin{thm}
\label{thm:A a F}A) The number $N$ of curves of genus $2$ on a
Fano surface of type $G$ is given by the following table: \smallskip{}
\\
\begin{tabular}{|c|c|c|c|c|c|c|c|}
\hline 
Group $G$ & $\mathbb{Z}/2\mathbb{Z}$ & $\mathbb{D}_{2}$ & $\mathbb{D}_{3}$ & $\mathbb{D}_{5}$ & $\mathbb{D}_{6}$ & $\mathbb{A}_{5}$ & $PSL_{2}(\mathbb{F}_{11})$\tabularnewline
\hline 
$N$  & 1 & 3 & 3 & 5 & 7 & 15 & 55\tabularnewline
\hline
\end{tabular}\\
For $S$ generic these genus $2$ curves are smooth.\\
B) The $3$ genus $2$ curves of a surface of type $\mathbb{D}{}_{2}$
are disjoint.\\
C) A surface $S$ of type $\mathbb{D}_{3}$ contains $3$ curves
$D_{1},D_{2},D_{3}$ of genus $2$, such that $D_{i}D_{j}=2$ if $i\not=j$.\\
 There exists a fibration $\gamma:S\rightarrow E$ onto an elliptic
curve such that $D_{1}+D_{2}+D_{3}$ is a fiber of $\gamma$.\\
D) A surface $S$ of type $\mathbb{D}_{5}$ contains $5$ curves
$D_{1},\dots,D_{5}$ of genus $2$, such that $D_{i}D_{j}=1$ if $i\not=j$.
\\
There exists a fibration $\gamma:S\rightarrow E$ onto an elliptic
curve  such that $D_{1}+\dots+D_{5}$ is a fiber of $\gamma$.\\
E) Let $S$ be a surface of type $\mathbb{D}_{6}$. There exists
a fibration $\gamma:S\rightarrow E$ onto an elliptic curve and genus
$2$ curves $D_{1},\dots,D_{7}$ such that $F_{1}=D_{1}+D_{2}+D_{3}$
and $F_{2}=D_{4}+D_{5}+D_{6}$ are $2$ fibers of $\gamma$ and $D_{7}$
is contained inside a third fiber. By C), the fibers $F_{1}$ and
$F_{2}$ corresponds to the $2$ subgroups $\mathbb{D}_{3}$ of $\mathbb{D}_{6}$.\\
F) A Fano surface $S$ of type $\mathbb{A}_{5}$ contains $15$
smooth curves of genus $2$. They generates a sub-lattice $\Lambda$
of $\NS(S)$ of rank $15$, signature $(1,14)$ and discriminant $2^{24}3^{6}$.
For $S$ generic $\Lambda$ has finite index inside $\NS(S)$. There
exist an infinite number of surfaces of type $\mathbb{A}_{5}$ with
maximal Picard number $h^{1,1}=25$.\end{thm}
\begin{rem}
In \cite{Roulleau1}, we give Fano surfaces of type of some groups
(by example the symmetric group $\Sigma_{5}$), but the genus $2$
curves we obtain in that way are sum of elliptic curves, and we are
interested by smooth genus 2 curves.
\end{rem}
The Klein cubic threefold: \[
F_{Kl}=\{x_{1}x_{5}^{2}+x_{5}x_{3}^{2}+x_{3}x_{4}^{2}+x_{4}x_{2}^{2}+x_{2}x_{1}^{2}=0\}\]
is the only one cubic of type $PSL_{2}(\mathbb{F}_{11})$ and this
group is its full automorphism group \cite{Adler}. It contains $55$
involutions. Let $S_{Kl}$ be the Fano surface of lines of $F_{Kl}$.
To each involution $g$ of $PSL_{2}(\mathbb{F}_{11})$, we denote
by $D_{g}\hookrightarrow S_{Kl}$ the corresponding curve of arithmetical
genus $2$ on $S$ (see Corollary \ref{cor:To-each-order}).
\begin{thm}
\label{thm:The-configuration-of}The $55$ genus $2$ curves $D_{g}$
are smooth. Their configuration is as follows : \begin{equation}
D_{g}D_{h}=\left\{ \begin{array}{ccc}
-4 & if & g=h\\
0 & if & o(gh)=2\,\, or\,\,6\\
2 & if & o(gh)=3\\
1 & if & o(gh)=5\end{array}\right.\label{eq:intersection Thm}\end{equation}
where we denote by $o(g)$ the order of an automorphism $g$.\\
 The sublattice $\Lambda'$ of the Néron-Severi group $\NS(S_{Kl})$
generated by the $55$ genus $2$ curves has rank $25=h^{1,1}(S_{Kl})$
and discriminant $2^{2}11^{10}$. \\
The group $\NS(S_{Kl})$ is generated by $\Lambda'$ and the class
of an incidence divisor.
\end{thm}
Let us prove Theorems \ref{thm:A a F} and \ref{thm:The-configuration-of}.

Let $D_{g},D_{h}$ be $2$ genus $2$ curves on a Fano surface $S$
corresponding (Corollary \ref{cor:To-each-order}) to involutions
$g$, $h$ and let $n$ be the order of $gh$.
\begin{lem}
\label{lem:The-intersection-number}Suppose that $n\in\{2,3,5,6\}$
and $S$ has type $\mathbb{D}_{n}$. The intersection number $D_{h}D_{g}$
is independent of the Fano surface of type $\mathbb{D}_{n}$.\end{lem}
\begin{proof}
The family $V$ of cubics forms $F_{eq}\in\mathbb{P}(H^{0}(\mathbb{P}^{4},\mathcal{O}(3)))$
such that the cubic $\{F_{eq}=0\}$ is smooth is open.\\
The group $\mathbb{D}_{n}$ acts naturally on $H^{0}(\mathbb{P}^{4},\mathcal{O}(3))$,
the third symmetric power of $H^{0}(\mathbb{P}^{4},\mathcal{O}(1))$.
For a cubic $F$ of type $\mathbb{D}_{n}$, there is a cubic form
$F_{eq}$ in $V$ such that $F=\{F_{eq}=0\}$ and there exists a character
$\chi:\mathbb{D}_{n}\rightarrow\mathbb{C}^{*}$ such that $F_{eq}\circ g=\chi(g)F_{eq}$.
For $n\in\{2,3,5,6\}$, the smoothness condition on $F$ implies that
$\chi$ is trivial (we used a computer). \\
Therefore, the family $V_{n}$ of cubic of type $\mathbb{D}_{n}$
is an open set of a projective space. Now it suffice to consider a
smooth curve inside $V_{n}$ between two points. That gives a flat
family of smooth cubic threefolds and of Fano surfaces. On each Fano
surfaces we get two genus $2$ curves and the two families of genus
$2$ curves are flat (see \cite{Hartshorne}, III, 9.7). The intersection
number of the genus 2 curves is therefore constant.
\end{proof}
Let $S$ be a Fano surface of type $\mathbb{D}_{5}$ and let $D_{g}$
and $D_{h}$ be as in the previous lemma.
\begin{lem}
\label{lem:intersection deux courbes eg 1}We have : $D_{g}D_{h}=1$.\end{lem}
\begin{proof}
Let $S_{0}$ be the Fano surface of the cubic $F_{0}=\{x_{1}^{3}+\dots+x_{5}^{3}=0\}$.
Let $\sigma$ be an element of the permutation group $\Sigma_{5}$
; $\sigma$ acts on $\mathbb{C}^{5}$ by $x\rightarrow(x_{\sigma1},\dots,x_{\sigma5})$
and acts on $F$ by taking the projectivisation. \\
 The involutions $g=(1,3)(4,5)$ and $h=(1,2)(3,5)$ have trace
1, their product is an order five element with trace equal to $0$.
In \cite{Roulleau}, we prove that the corresponding genus $2$ curves
are $D=E_{1}+E_{2}$ and $D'=E'_{1}+E'_{2}$ with $E{}_{1},E{}_{2},E'_{1},E'{}_{2}$
elliptic curves on $S_{0}$ such that $DD'=1$. We now apply Lemma
\ref{lem:The-intersection-number} to $S$ of type $\mathbb{D}_{5}$.
\end{proof}
Let $S_{Kl}$ be the Fano surface of the Klein cubic. Let $x,y,z,w$
be elements of $\{0,1,2\}$ and let $\Lambda_{x,y,z,w}$ be the lattice
generated by $55$ generators $L_{g}$ with intersection numbers :\[
L_{g}L_{h}=\left\{ \begin{array}{ccc}
-4 & if & g=h\\
x & if & o(gh)=2\\
y & if & o(gh)=3\\
z & if & o(gh)=5\\
w & if & o(gh)=6\end{array}\right..\]
By Theorem \ref{thm:A-curve-D} and Lemma \ref{lem:The-intersection-number},
the lattice generated by the $55$ genus $2$ curves on $S_{Kl}$
is isomorphic to one of the lattices $\Lambda_{x,y,z,w}$. The lattices
$\Lambda_{0,2,1,0}$ and $\Lambda_{0,0,0,2}$ are the only ones of
rank less or equal to $25=h^{1,1}(S)$. By Lemma \ref{lem:intersection deux courbes eg 1},
the lattice $\Lambda_{0,0,0,2}$ cannot be the lattice generated by
the $55$ genus $2$ curves on $S_{Kl}$.\\
The lattice $\Lambda_{0,2,1,0}$ has rank $25$, discriminant $2^{2}11^{10}$,
signature $(1,24)$. In \cite{Roulleau2}, we computed a basis of
$\NS(S_{Kl})$ : it has rank $25$ and discriminant $11^{10}$. The
lattice $\Lambda_{0,1,2,0}$ has thus index $2$ in $\NS(S_{Kl})$
and, as we know that $C_{s}D_{g}=2$, we can check that the incidence
divisor $C_{s}$ is not in this lattice. 

By Lemma \ref{lem:The-intersection-number}:
\begin{cor}
\label{cor:The-formula-intersection}The formula \ref{eq:intersection Thm}
giving the intersection number of genus $2$ curves on $S_{Kl}$ holds
also for the Fano surfaces of type $\mathbb{D}_{n},n\in\{2,3,5,6\}$
and $\mathbb{A}_{5}$.\end{cor}
\begin{proof}
The groups $\mathbb{D}_{n},n\in\{2,3,5,6\}$ are subgroups of $\mathbb{A}_{5}$
and $PSL_{2}(\mathbb{F}_{11})$.
\end{proof}
Let $S$ be a Fano surface of type $\mathbb{D}_{3}$. Let us prove
that:
\begin{lem}
There exists a fibration $\gamma:S\rightarrow E$ with $E$ an elliptic
curve such that $D_{1}+D_{2}+D_{3}$ is a fiber of $\gamma$. \end{lem}
\begin{proof}
We have : $(D_{1}+D_{2}+D_{3})^{2}=0$. \\
The dihedral group $\mathbb{D}_{3}$ acts on $H^{0}(\Omega_{S})$
by two copies of a representation of degree $2$ plus the trivial
representation of degree $1$ (see also the next paragraph). \\
Let $g$ be an involution of $\mathbb{D}_{3}$. We know its action
on $H^{0}(\Omega_{S})^{*}$, in particular, we know the eigenspace
$T(g)$ with eigenvalues $-1$. By Lemma \ref{lem: espace tangent ourbed},
$T(g)$ is the tangent space of a genus $2$ curve. We can check that
the sub-space $W$ of $H^{0}(\Omega_{S})^{*}$ generated by the tangent
space coming from the $3$ genus $2$ curves has dimension $4$. \\
The space $W$ is the tangent space of a $4$ dimensional variety
$B$ inside the Albanese variety of $S$. Thus there exist a fibration
$q:\Alb(S)\rightarrow E$ where $E$ is an elliptic curve. The $3$
genus $2$ curves are contracted by the composition of the Albanese
map and $q$.
\end{proof}
The similar assertions D) and E) of Theorem \ref{thm:A a F} relative
to the fibrations of surfaces of type $\mathbb{D}_{5}$ and $\mathbb{D}_{6}$
are proved in the same way. 
\begin{lem}
For $S$ generic of type $G=\mathbb{Z}/2\mathbb{Z}$ or $\mathbb{D}_{n},\, n\in\{2,3,5,6\}$,
the genus $2$ curves $D_{g}$ are smooth.\end{lem}
\begin{proof}
It is enough to check that the conic $\mathcal{Q}$ defined in paragraph
\ref{sub:Involutive-automorphisms-and} is smooth if the cubic of
type $G$ is generic (we give the equation of such cubic in the next
paragraph). 
\end{proof}
Let us prove part F) of Theorem \ref{thm:A a F}. Let $S$ be a Fano
surface of type $\mathbb{A}_{5}$. By Corollary \ref{cor:The-formula-intersection},
we know the intersection numbers of the $15$ genus $2$ curves $D_{g}$
corresponding to the $15$ involutions of $\mathbb{A}_{5}$. The lattice
$\Lambda$ generated by these $15$ curves has signature $(1,14)$
and discriminant $2^{24}3^{6}$.
\begin{lem}
The genus $2$ curves $D_{g}$ on a Fano surface $S$ of type $\mathbb{A}_{5}$
and $PSL_{2}(\mathbb{F}_{11})$ are smooth.\end{lem}
\begin{proof}
In \cite{Roulleau}, we proved that the elliptic curves $E$ on a
Fano surface $S$ correspond bijectively with automorphisms $\sigma_{E}\in\Aut(S)$
such that the trace of $\sigma_{E}$ acting on $H^{0}(\Omega)$ is
$-3$ (involutions of type I). We classified all automorphism groups
generated by involutions of type I.\\
Moreover, by \cite{Roulleau}, Theorem 13, for $2$ elliptic curves
$E\not=E'$ on $S$, we have $(\sigma_{E}\sigma_{E'})^{2}=1$ if and
only if $EE'=1$ i.e. if and only if $E+E'$ is a genus $2$ curve
on $S$. In that case, the trace of the involution $\sigma_{E}\sigma_{E'}$
on $H^{0}(\Omega)$ is $1$ (involution of type II). \\
Suppose that one genus $2$ curve on a Fano surface of type $\mathbb{A}_{5}$
is the sum of $2$ elliptic curves. By transitivity, the $15$ genus
$2$ curves are also sum of genus $2$ curves and $\mathbb{A}_{5}$
(group generated by automorphisms of type II), is an automorphism
sub-group of the group generated by involutions of type I.\\
 By the classification of automorphism groups generated by involutions
of type I (\cite{Roulleau}, Theorem 26), $\mathbb{A}_{5}$ is a subgroup
of the symmetric group $\Sigma_{5}$ or of the reflection group $G(3,3,5)$
acting on $\mathbb{C}^{5}$. Thus it must exists elements $a,b$ of
$\Sigma_{5}$ or $G(3,3,5)$ such that: \[
a^{2}=b^{3}=(ab)^{5}=1\]
(relations defining the alternating group of degree $5$) and $Tr(a)=1$,
$Tr(b)=-1$. But it is easy to check that no such elements exist in
$\Sigma_{5}$ nor in $G(3,3,5)$. Therefore the $15$ genus $2$ curves
on a Fano surface of type $\mathbb{A}_{5}$ are smooth.\\
As an involution in $PSL_{2}(\mathbb{F}_{11})$ in contained inside
a group $\mathbb{A}_{5}$, the $55$ genus $2$ curves on $S_{Kl}$
are smooth.
\end{proof}
Let $a,b$ be the generators $a,b$ of $\mathbb{A}_{5}$ such that:
\[
a^{2}=b^{3}=(ab)^{5}=1,\]
with $a$ the diagonal matrix with diagonal elements $-1,-1,1,1,1$
and:\[
b=\left(\begin{array}{ccccc}
-\frac{1}{2} & -\frac{1}{2} & \frac{1}{2} & -\frac{1}{2} & 0\\
\frac{1}{2} & 0 & \frac{1}{2} & 0 & -\frac{1}{2}\\
-\frac{1}{2} & \frac{1}{2} & \frac{1}{2} & \frac{1}{2} & 0\\
\frac{1}{2} & 0 & \frac{1}{2} & 0 & \frac{1}{2}\\
0 & 0 & -2 & -2 & -1\end{array}\right).\]
A cubic threefold of type $\mathbb{A}{}_{5}$ has equation:\[
\begin{array}{c}
a(x_{4}^{3}+x_{4}(x_{1}^{2}-x_{2}^{2}+x_{3}^{2})+x_{3}(-x_{2}^{2}+3x_{4}^{2}+x_{5}^{2})+2x_{3}^{2}x_{5}\\
+2x_{1}x_{2}(x_{3}+x_{4}+x_{5})+2x_{3}x_{4}x_{5})+b(-x_{3}^{3}+x_{3}(x_{1}^{2}-x_{2}^{2}-x_{4}^{2})\\
+x_{4}(x_{1}^{2}-3x_{3}^{2}-x_{5}^{2})-2x_{4}^{2}x_{5}+2x_{1}x_{2}(x_{3}+x_{4}+x_{5})-2x_{3}x_{4}x_{5})=0.\end{array}\]
Let $t$ be the point on $S$ corresponding to the line $\{x_{3}=x_{4}=x_{5}=0\}$.
The Quintic $\Gamma_{t}$ is $\Gamma_{t}=Q+C$ for $C$ the cubic:\[
\begin{array}{c}
a(x_{4}^{3}+x_{4}x_{3}^{2}+x_{3}(3x_{4}^{2}+x_{5}^{2})+2x_{3}^{2}x_{5}+2x_{3}x_{4}x_{5})\\
+b(-x_{3}^{3}-x_{3}x_{4}^{2}-x_{4}(3x_{3}^{2}+x_{5}^{2})-2x_{4}^{2}x_{5}-2x_{3}x_{4}x_{5})=0.\end{array}\]
This define a pencil of elliptic curves and by Remark \ref{rem: sur les theta caracteristiques},
the cubics of type $\mathbb{A}{}_{5}$ form a $1$ dimensional family
of cubics.\\
Let $\Alb(S)$ be the Albanese variety of $S$. By the symmetries
of $\mathbb{A}_{5}$ acting on $\Alb(S)$, this Abelian variety is
isogenous to $E^{5}$ for some elliptic curve $E$. If $E$ has no
complex multiplication, then the Picard number of $\Alb(S)$ is $15$,
otherwise it is $25$. By \cite{Roulleau3}, the Picard number of
a Fano surface $S$ and of its Albanese variety $\Alb(S)$ are equal.
Thus the Picard number of $S$ is $15$ or $25$. As the curve $E$
varies, the case $E$ with CM occurs.
\begin{rem}
The number of elliptic curves on a Fano surface is bounded by $30$
and the Fano surface of the Fermat cubic threefold is the only one
to contain $30$ elliptic curves. It is tempting to think that $55$
is the bound for the number of smooth genus $2$ curves on a Fano
surface and that $S_{Kl}$ is the only one to reach this bound.
\end{rem}

\subsection{Construction of Fano surfaces of a given type, on the completeness
of the groups classification.}

In order to get a classification of automorphism groups of Fano surfaces
generated by involutions $g,h,\dots$, it is natural to study their
products i.e. to look at the dihedral group generated by two elements
$g,h$. In \cite{Roulleau2}, we prove that the order of $Aut(S)$
is prime to $7$ and that an automorphism $\sigma$ of $S$ preserve
a $5$ dimensional principally polarized Abelian variety. That implies
(\cite{Birkenhake}, Proposition 13.2.5 and Theorem 13.2.8) that the
Euler number of the order $n$ of $\sigma$ is less or equal to $10$,
therefore: \[
n\in\{2,\dots,6,8,9,10,11,12,15,16,18,20,22,24,30\}.\]
Hence, it is wise to study representations of dihedral group of order
$2n$ with $n$ in the above set, such that the order $2$ elements
have trace equal to $1$. Our method is a case by case check ; we
have results for the cubic threefold of type $\mathbb{D}_{n}$ with
$n\in\{2,3,4,5,6,8,11,12,16,20,22,24\}$.

The $V_{\frac{k}{n}},\,0<k<n$ representation of the dihedral group
of order $2n$ (generated by $a,b$ such that $a^{n}=b^{2}=1,\, bab=a^{-1}$)
is given by the matrices:\[
a=\left(\begin{array}{cc}
\cos\frac{2k\pi}{n} & -\sin\frac{2k\pi}{n}\\
\sin\frac{2k\pi}{n} & \cos\frac{2k\pi}{n}\end{array}\right),\, b=\left(\begin{array}{cc}
1 & 0\\
0 & -1\end{array}\right).\]
There is also the trivial representation $T$, the linear representation
$ $$L:\, a\rightarrow1,b\rightarrow-1$, and if $n$ is even, the
representations $L_{1}:\, a\rightarrow-1,b\rightarrow1$ and $L_{2}:\, a\rightarrow-1,b\rightarrow-1$.
The representation $V_{\frac{k}{n}}$ is faithful if and only if $k$
is prime to $n$ ; it is irreducible if and only if $k\not=\frac{n}{2}$.
The representations $V_{\frac{k}{n}}$ and $V_{\frac{n-k}{n}}$ are
equivalent.

$\circledcirc$ A surface of type $\mathbb{Z}/2\mathbb{Z}$ is the
Fano surface of the cubic threefold $F_{2}$ given in paragraph \ref{sub:Involutive-automorphisms-and}.

$\circledcirc$ The group $\mathbb{D}_{2}$ is given by the representation
$L+L_{1}+L_{2}+2T$. An invariant cubic $F_{4}$ has equation:\[
F_{4}=\{(ax_{1}^{2}+bx_{2}^{2}+cx_{3}^{2})x_{4}+(dx_{1}^{2}+ex_{2}^{2}+fx_{3}^{2})x_{5}+gx_{4}^{3}+hx_{5}^{3}+kx_{1}x_{2}x_{3}=0\}\]
where $a,\dots,k$ are constants. 

$\circledcirc$ The representation $\mathbb{D}_{3}$ is given by $2V\frac{1}{3}+T$.
A cubic threefold of type $\mathbb{D}{}_{3}$ has equation:\[
\begin{array}{c}
x_{5}^{3}+(x_{1}^{2}x_{5}+x_{2}^{2}x_{5})+(x_{3}^{2}x_{5}+x_{4}^{2}x_{5})+a(x_{1}^{3}-3x_{2}^{2}x_{1})+b(x_{3}^{3}-3x_{4}^{2}x_{3})+c(x_{1}x_{3}x_{5}\\
+x_{2}x_{4}x_{5})+d(x_{3}^{2}x_{1}-x_{4}^{2}x_{1}-2x_{2}x_{3}x_{4})+e(x_{1}^{2}x_{3}-x_{2}^{2}x_{3}-2x_{1}x_{1}x_{4}).\end{array}\]
Let us denote by $\mathbb{D}'_{3}$ the dihedral group of order $6$
such that an element of order $2$ (resp. $3$) has trace equals to
$1$ (resp $2$). Its representation is $V_{\frac{1}{2}}+L+2U$. A
cubic threefold $F$ of type $\mathbb{D}'_{3}$ has equation:\[
ax_{4}^{3}+bx_{5}^{3}+(x_{1}^{2}+x_{2}^{2})(ux_{4}+vx_{5})+c(x_{1}^{3}-3x_{2}^{2}x_{1})+dx_{5}^{2}x_{4}+ex_{4}^{2}x_{5}+fx_{3}^{2}x_{5}+gx_{3}^{2}x_{4}=0\]
The involutions $x\rightarrow(x_{1},\pm x_{2},\pm x_{3},x_{4},x_{5})$
act on $F$. This gives a genus $2$ curve on the Fano surface that
is sum of two elliptic curves : this is not interesting.\\
There is a third representation of the dihedral group of order
$6$ such that the trace of the order $2$ elements is $1$, but it
is not faithfull.

$\circledcirc$ The representation $\mathbb{D}_{5}$ is given by $V\frac{1}{5}+V\frac{2}{5}+T$.
A cubic threefold of type $\mathbb{D}{}_{5}$ has equation:\[
\begin{array}{c}
x_{5}^{3}+c(x_{1}^{2}x_{5}+x_{2}^{2}x_{5})+d(x_{3}^{2}x_{5}+x_{4}^{2}x_{5})+a(x_{1}^{2}x_{3}-x_{2}^{2}x_{3}\\
+2x_{1}x_{2}x_{4})+b(-x_{3}^{2}x_{1}+x_{4}^{2}x_{1}+2x_{2}x_{3}x_{4})=0.\end{array}\]

$\circledcirc$ The representation $\mathbb{D}_{6}$ is given by $V\frac{1}{6}+V\frac{2}{6}+T$.
A cubic threefold of type $\mathbb{D}{}_{6}$ has equation:\[
ax_{5}^{3}+b(x_{1}^{2}x_{5}+x_{2}^{2}x_{5})+c(x_{3}^{2}x_{5}+x_{4}^{2}x_{5})+d(x_{3}^{3}-3x_{4}^{2}x_{3})+e(x_{1}^{2}x_{3}-x_{2}^{2}x_{3}+2x_{1}x_{2}x_{4})=0.\]
The dihedral group of order $12$ contain the dihedral group of order
$6$. For this last group, the only interesting representation is
$\mathbb{D}_{3}$ (such that the trace of an order $3$ element equals
$-1$). With that point in mind, we can check that the only interesting
representation of the dihedral group of order $12$ is $\mathbb{D}_{6}$. 

$\circledcirc$ There is another $5$ dimensional representation of
the alternating group of degree $5$ such that the trace of an order
$2$ element is $1$, but the trace of an order $3$ element is $2$
as for $\mathbb{D}'_{3}$ and that implies that the corresponding
genus $2$ are sum of $2$ elliptic curves : this is not interesting.

$\circledcirc$ Let us now prove the following
\begin{prop}
There do not exist a Fano surface of type the dihedral group of order
$2n$ with $n\in\{4,8,12,16,20,24\}$.\end{prop}
\begin{proof}
Let $a,b$ be generators of the dihedral group of order $8$ such
that $a^{4}=1,b^{2}=1,bab=a^{3}$. We are looking for representations
such that the trace of the order $2$ elements $b,ab,a^{2}b,a^{3}b$
is $1$. For the traces of $a$ and $a^{2}$, the possibilities are
$Tr(a)=-1,Tr(a^{2})=1$ or $Tr(a)=3,Tr(a^{2})=1$ or $Tr(a)=1,Tr(a^{2})=-3$.
In each cases, we computed the spaces $V_{\chi}$ of cubics such that
$F_{eq}\circ g=\chi(g)F_{eq}$ for character $\chi$. That gives no
smooth cubic threefolds.\\
The dihedral groups of order $16,24,32,40,48$ contain the dihedral
group of order $8$, thus they cannot be type of a cubic.
\end{proof}

\subsection{A conjecture.\protect \\
}

The lattice $\Lambda_{0,0,0,2}$ of the proof of Theorem \ref{thm:The-configuration-of}
has rank $21$, discriminant ${11.2}^{22}$ and signature $(1,20)$.
It is remarkable that this lattice has the right signature to be the
Néron-Severi group of a surface. \\
There is a natural representation of the group $PSL_{2}(\mathbb{F}_{11})$
on $\mathbb{C}^{5}$ (for which the Klein cubic generates the space
invariant cubics). On the other hand classification of surfaces not
of general type in $\mathbb{P}^{4}=\mathbb{P}(\mathbb{C}^{5})$ is
an old unsolved problem.\\
The author's opinion is that we should consider the lattice $\frac{1}{2}\Lambda_{0,0,2,0}$
of rank $21$ given by the generators $\{L_{g},g$ involution$\}$
and the relations: \[
L_{g}L_{h}=\left\{ \begin{array}{ccc}
-2 & if & g=h\\
1 & if & o(gh)=6\\
0 &  & otherwise\end{array}\right.\]
 as the lattice generated by a configuration of $55$ lines on a surface
$S^{?}$ in $\mathbb{P}^{4}$. That surface $S^{?}\hookrightarrow\mathbb{P}^{4}$
is expected as the (may be non-complete) intersection of invariants
of the group $PSL_{2}(\mathbb{F}_{11})$ acting on $\mathbb{P}^{4}$,
in such a way that $PSL_{2}(\mathbb{F}_{11})$ acts on it. \\
On a surface of general type, the lattice generated by (-2)-curves
is negative definite. As the lattice $\frac{1}{2}\Lambda_{0,0,2,0}$
is generated by (-2)-curves and has signature $(1,20)$, the surface
$S^{?}$ can not be of general type.

Xavier Roulleau\\
Graduate School of Mathematical Sciences, University of Tokyo,
3-8-1 Komaba, Meguro, Tokyo, 153-8914 Japan.\\
roulleau@ms.u-tokyo.ac.jp
\end{document}